\documentclass[11pt]{article}
\usepackage{geometry}                
\usepackage{cyrillic}                  
\usepackage{graphicx}
\usepackage{amssymb}
\usepackage{epstopdf}
\usepackage{amssymb,latexsym,graphicx,psfrag,epstopdf,hyperref}
\usepackage[all,arc,dvips]{xy}
\newcommand{\cyrrm}{\fontencoding{OT2}\selectfont\textcyrup}

\usepackage{amsthm}

\newtheorem{mydef}{Definition}
\newtheorem{mythm}{Theorem}
\newtheorem{myfact}{Fundamental Fact}

\title{Volumes in Hyperbolic Space}
\author{Christina Laternser}
\date{Spring 2010}                                           

\begin{document}

\maketitle
\begin{abstract}
This paper focuses on the investigation of volumes of large Coxeter hyperbolic polyhedron. First, the paper investigates the smallest possible volume for a large Coxeter hyperbolic polyhedron and then looks at the volume of pyramids with one vertex at infinity. 
\end{abstract}
 
\section {Introduction}
 \indent This research started in the summer of 2009 with a fellow student, Michael Fry who had graduated that summer. During those summer months, we spent our time investigating volumes of large Coxeter hyperbolic polyhedron. The main focus was to explore hyperbolic geometry beyond the class room and to see if we could find the smallest possible volume for this class of objects. By the end of the summer, we had found a polyhedron which generated the smallest volume we had seen. The polyhedron we had found was the Lambert cube. Following this research, the fall of 2009 and spring of 2010, were spent on further research by myself. The goal for both semesters was supposed to lead into my being introduced to Hyperbolic Dehn Filling of ideal polyhedra. We did not get to spend as much time on it as needed, and therefore that part of my research will not be included in this research paper. \\
\indent To have a better and more formal understanding of Hyperbolic Geometry, we were introduced to Curved Spaces, a program which visualizes Hyperbolic Geometry, and Orb, a program which calculates volumes of hyperbolic polyhedra. On a side quest, I tried to combine both concepts of Curved Spaces and Orb. I wanted a program which would allow for the visualization of  any kind of hyperbolic polyhedra given the matrices which were generated by Orb. Curved Spaces only had five preprogrammed hyperbolic polyhedra, and to see the polyhedra one is working with makes working with the polyhedra a lot easier. However, I did not have a chance to finish the program, which is why that part is also not included in this research paper. 
\section{Definitions}

This paper assumes that the reader is familiar and comfortable with Euclidean Geometry. Some concepts therefore will be assumed to be known to the reader and most likely not explicitly stated. Euclidean Geometry, as the name suggests, was named for Euclid (300 B.C.). Euclid established five postulates as the foundation of Euclidean Geometry:
\\
\begin{enumerate}
\item Given any two points in space, there exists a unique line segment connecting these two points.

\item Any line segment can be extended indefinitely. 

\item Given any point in space, a circle can be drawn with any radius around that point. 

\item Any two right angles are congruent.

\item (Parallel Postulate) Given any point $P$ in space and a line $l_1$ that does not pass through $P$, there exists a unique line $l_2$ through $P$ which does not intersect $l_1$. \cite{ab12}
\\
\end{enumerate}

\indent Mathematics is meant to describe the world we live in. Euclidean Geometry, or what people most commonly refer to as just Geometry, describes the physical world in fundamental ways. Euclid observed these things and decided that these few postulates are all that is needed in order to successfully describe everything around us. We tend to assume that we live in a Euclidean Geometry. Notions of distance and concepts of triangles are well defined to and by us. We can find the shortest distance between two points by using the Pythagorean Theorem, assuming, of course, that we superimpose a grid on the plane the two points lie in. A triangle is simply defined to be three points in space which are joined by three line segments. We also know that the sum of the interior angles of a triangles sum up to $\pi$. There are also fundamental operations called isometries which are the motions of the plane that preserve distances and angles, and whose existence is added to the axioms of Euclidean geometry. There are three isometries: translations, rotations, and reflections. \cite{ab12} A translation can be thought of as the shifting of the plane by a scalar; translations have no fixed points or is what is called the identity. A rotation fixes one point $P$ in space and rotates space around $P$ by a $\theta$ degrees. A reflection fixes a line $l$ in space and reflects every point through $l$. The only points that are not moved by the reflection are the points that are on $l$. \\
\indent When people are first introduced to Non-Euclidean Geometry, it comes as a shock, because people are rarely taught anything outside of Euclidean Geometry. A question to ask yourself is, do we actually live in a world that follows the rules Euclid set down over two thousand years ago, or do we live in another kind of geometry? It is hard to imagine that we live in anything but Euclidean Geometry. At the same time, we are very small and the world is very big. For a very long time, we did not even know that the earth was round! Mathematicians found ways to describe the world we live in, found geometries which followed all of Euclid's postulates, except the 5th one; and so, Hyperbolic Geometry was born. In Hyperbolic Geometry, all of Euclid's postulates hold, except the Parallel Postulate. The revised 5th postulate for Hyperbolic Geometry goes as follows: \textquotedblleft Given any point $P$ in space and a line $l_1$, there are { \em infinitely} many lines through $P$ which are parallel to $l_1$'' \cite{ab12}.
\\
\\
\indent Envisioning the hyperbolic plane, $\mathbb{H}^2$, is for the most part impossible, hence models need to be used in order to work with $\mathbb{H}^2$ or any higher dimensions. In the case of $\mathbb{H}^2$ the Poincar\'e disk (see Figure 3 for an illustration) or the Half Plane model (see Figure 2 for an illustration) are often used and this paper will refer to either the Poincar\'e disk/ sphere or the Half Plane/Space model in order to illustrate images. 
 
Due to the revised Parallel Postulate, there are certain properties which no longer hold and other properties which are completely new and are unique to Hyperbolic Geometry:
\\

\begin{enumerate}
\item A scaling of $\mathbb{E}$ is an operation that multiplies all distances by a constant factor. Scaling objects no longer exist in Hyperbolic Geometry. Because of the Revised Parallel Postulate, changing the length of an edge results in an entirely new figure. 
\item Instead of only the rectangles, 30-60-90 triangles, 45-45-90 triangles, and regular hexagons tiling the Euclidean plane, there are infinitely many distinct polygons which tile the (Hyperbolic) plane. For example, and square with interior angles of $\frac{\pi}{4}$ tiles out $\mathbb{H}^2$. The tiling of the plane (or space), also referred to as tessellation, is a collection of figures that fill space without overlapping or having gaps. A polygon that tiles out the plane is referred to as a Coxeter polygon. 
\item The sum of the interior angles of a triangle in Hyperbolic Geometry no longer has to sum up to $\pi$. In fact, only occasionally will the sum of the interior angles of a triangle add up to $\pi$; we call these triply asymptotic triangles. In all other cases, the sum will be less than $\pi$. More specifically, the area of any triangle in Hyperbolic space is the difference between $\pi$ and the interior angles. 
\item Lines are paths that minimize distance between points on them. So, for example, in the upper half plane model, lines of $\mathbb{H}^2$ come in two varieties, vertical Euclidean lines and arcs of semicircles perpendicular to the x-axis (see Figure 1). 
\\
\end{enumerate}

\begin{figure}[htbp] 
   \centering
   \includegraphics[width=4in]{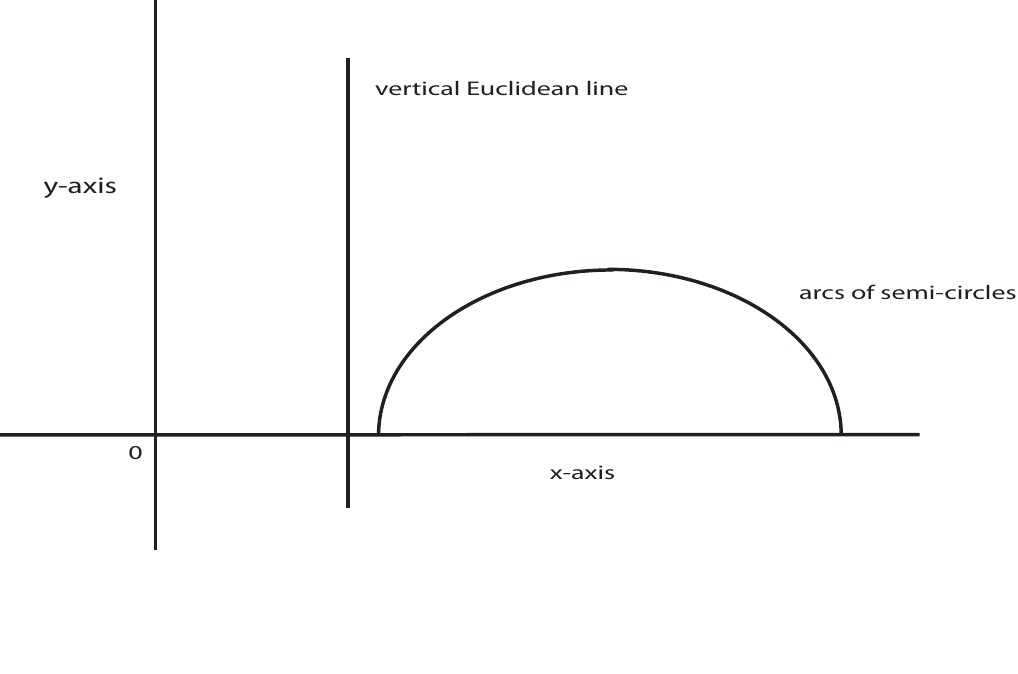} 
   \caption{Lines in the Half Space model }
   \label{fig:example}
\end{figure}

\indent These are just a few examples of things which change when working in Hyperbolic Geometry. The lines in the upper half-plane model allow us to easily visualize the need to modify Euclid's 5th postulate. 
\begin{figure}[htbp] 
   \centering
   \includegraphics[width=4in]{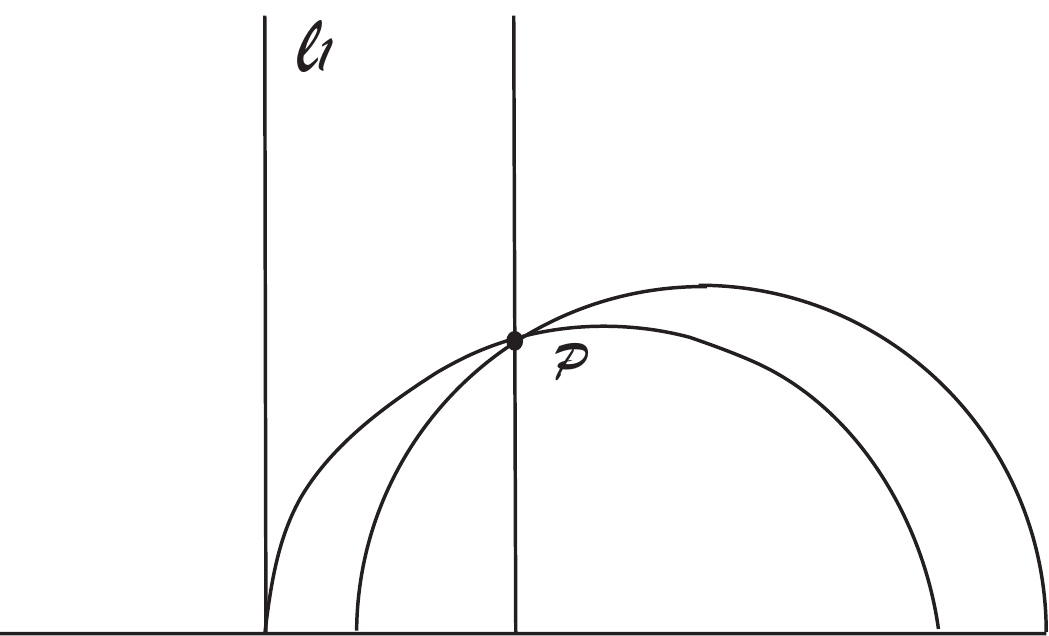} 
   \caption{Revised Parallel Postulate in the Half-plane model }
   \label{fig:example}
\end{figure}

\indent In $\mathbb{E}^2$, scaling or resizing a polygon is possible. For example, a square with side length 1 tiles out the plane. But a square with side length $\frac{1}{2}$ also tiles out the plane, so does a square of side length $\frac{1}{4}$, and so forth. In Hyperbolic space, scaling a square is not possible. Changing the length of an edge may change the polygon into a polygon that no longer tiles out $\mathbb{H}^2$. Every polygon that is found to tile out $\mathbb{H}^2$ is unique, because there is a specific angle associated with each vertex, and changing the angle between the edges creates a different polygon. The same concept works for $\mathbb{H}^3$. 

\begin{figure}[htbp] 
   \centering
   \includegraphics[width=4in]{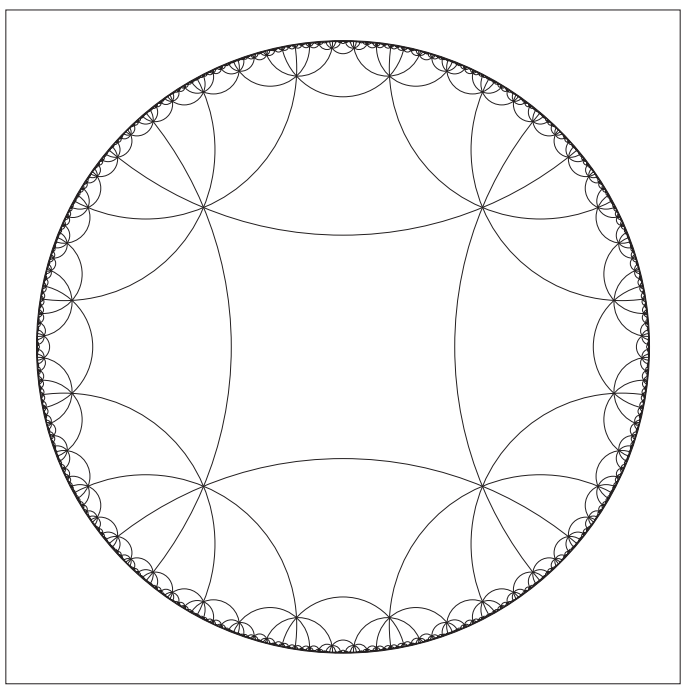} 
   \centering
   \caption{This image was generated using the program \emph{lim}, which was created by Curt McMullen. It depicts the tiling of $\mathbb{H}^2$ by a square with interior angles $\frac{\pi}{4}$.}
   \label{fig:example}
\end{figure}

\indent It is convenient to work with the planar graphs obtained by projecting the 1-dimensional skeleton of our polyhedra into the plane stereographically. An abstract polyhedron is a trivalent graph in the plane which is made up of the edges of the polyhedron together with the regions of the plane bounded by the graph (called the faces). The unbounded region is also a face. An Abstract Polyhedron must satisfy these four conditions. 

\begin{mydef} \emph{\textbf{(Abstract Polyhedron)}}
\begin{enumerate}
\item If two faces intersect, then they do so in either one edge or one vertex. 
\item There are no holes, i.e. the abstract polyhedra is closed.
\item It must have more than three faces.
\item Each edge meets exactly two faces.\cite{CHK}
\end{enumerate} 
\end{mydef}
It follows that in each vertex exactly three edges meet. 
Before I can talk about the most important theorem in my research, I must define a couple of terms, which Andreev's theorem mentions. \\
\indent Firstly, the angle between two planes is a dihedral angle. The numbers which the edge labels are labeled with are the dihedral angles the two faces meet in. When labeling the edges of an abstract polyhedron, it is convention to enumerate the edge by the denominator of the fraction. For example, instead of labeling an edge with $\frac{\pi}{2}$, it will be labeled simply by 2. Differently labelings of the abstract polyhedron will be discussed in more detail below. \\
A simple closed curve $\Gamma$ that crosses $k$ edges of an abstract polyhedron $C$ is called a $\emph{k-circuit}$. If all of the endpoints of the edges of $C$ which $\Gamma$ crosses are distinct then the circuit is considered to be a $\emph{prismatic k-circuit}$. A subset $S$ of $\mathbb{H}^3$ is $\emph{compact}$ if it is closed and bounded. \\
\indent Although many different kinds of polyhedron tile $\mathbb{H}^3$, not all do. Andreev's Theorem helps us identify which ones tile out $\mathbb{H}^3$ and which ones do not.

\begin{mythm} \emph{\textbf{(Andreev's Theorem)}}: 
Let $C$ be an abstract polyhedron with more than four faces and suppose that non-obtuse angles $\alpha_{i}$ are given corresponding to each edge $e_{i}$ of $C$. There is a compact hyperbolic polyhedron $P$ whose faces realize $C$ with dihedral angle $\alpha_{i}$ at each edge $e_{i}$ if and only if the following five conditions all hold:  
\begin{enumerate}
\item For each edge, $e_i$, the corresponding angle, $\alpha_{i}$, has to be larger than 0. 
\item Whenever three distinct edges, call them $e_i, e_j, e_k$, meet in a vertex, their angles, $\alpha_{i}$, $\alpha_{j}$, $\alpha_{k}$ must add up to greater than $\pi$. 
\item Whenever $\Gamma$ is a prismatic 3-circuit intersection edges  $e_i, e_j, e_k$, then the sum of the angles, $\alpha_{i}$, $\alpha_{j}$, $\alpha_{k}$, must add up to less than $\pi$.
\item  Whenever $\Gamma$ is a prismatic 4-circuit intersection edges  $e_i, e_j, e_k, e_l$, then the sum of the angles, $\alpha_{i}$, $\alpha_{j}$, $\alpha_{k}$, $\alpha_{l}$ must add up to less than 2$\pi$.
\item Whenever there is a four sided face bounded by edges $e_1,e_2,e_3$, and $e_4$ enumerated successfully, with edges $e_{12}, e_{23}, e_{34}, e_{41}$ entering the four vertices (edge $e_{ij}$ connects to the ends of $ e_i$ and $e_j$), then: 
\\
\\
\centering {$\alpha_1 + \alpha_3 + \alpha_{12} + \alpha_{23} + \alpha_{34} + \alpha_{41} < 3\pi$ , 
and

$\alpha_2 + \alpha_4 + \alpha_{12} + \alpha_{23} + \alpha_{34} + \alpha_{41} < 3\pi$} \footnote{ The last part of Andreev's Theorem only applies to Triangular Prisms. Hence, if the abstract polyhedron $P$ has more than five faces, then you must not concern oneself with the last part of Andreev's Theorem.}

\end{enumerate}
Furthermore, this polyhedron is unique up to isometries of $\mathbb{H}^3$.  \cite{rr12}
\end{mythm}

If a polyhedra satisfies all of Andreev's conditions, then the polyhedra is said to tile out $\mathbb{H}^3$. Proof of this theorem can be found in Roeder's paper on Andreev's theorem.$\emph{\cite{rr12}}$
\\
\\
Any time three edges meet they must satisfy one of the following in order to meet Andreev's 2nd, 3rd or 4th condition. 
\begin{center}
\begin{tabular}{|l | r | r | r|}
\hline
(2,3,3)	& tetrahedral 		\\ \hline
(2,3,4)	& cubal, octahedral 	\\ \hline
(2,3,5) 	& dodecahedral	\\ \hline 
(2,2,n) 	& dihedral			\\ \hline
\end{tabular}
\end{center}
\begin{figure}[htbp] 
   \centering
   \includegraphics[width=2.5in]{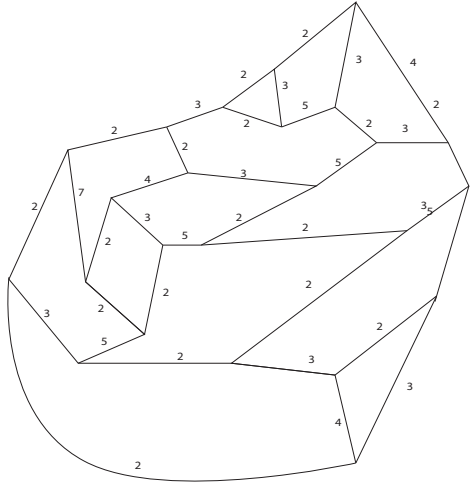} 
   \caption{Planar Graph with appropriate edge labels}
   \label{fig:example}
\end{figure}

Edge labelings which do not follow the restrictions given above will not generate a hyperbolic  Coxeter polyhedron. 

\begin{mydef}
A $\emph{topological disk}$ is a 2-dimensional region that can be contracted into a point. 
\end{mydef}

For our purposes, it suffices to think of a 2-orbifold as simple closed curve on a (abstract) polyhedron. 

\begin {mydef}
A $\emph{2-orbifold}$ in an abstract polyhedron is a simple closed curve together with the topological disk that it bounds. 
\end{mydef}

\begin{mydef} Let $F = D \cup C$ be a 2-orbifold in an abstract polyhedron, where $C$ is a simple closed curve and $D$ is the topological disk it bounds. A compression of $F$ is an arc $\alpha$ on the topological disk $D$ with end points on $C$, a 2-orbifold, that crosses 1 or 0 edges of the abstract polyhedron and such that the two arcs of $C$ on either side of $\alpha$ each cross the edges of the polyhedron in at least two places defines a compression of a 2-orbifold. A 2-orbifold is compressible if either it takes one of the forms in Figure 4 or it has a compression. Otherwise $F$ is called incompressible. 
\end{mydef}

\begin{figure}[htbp] 
   \centering
   \includegraphics[width=4in]{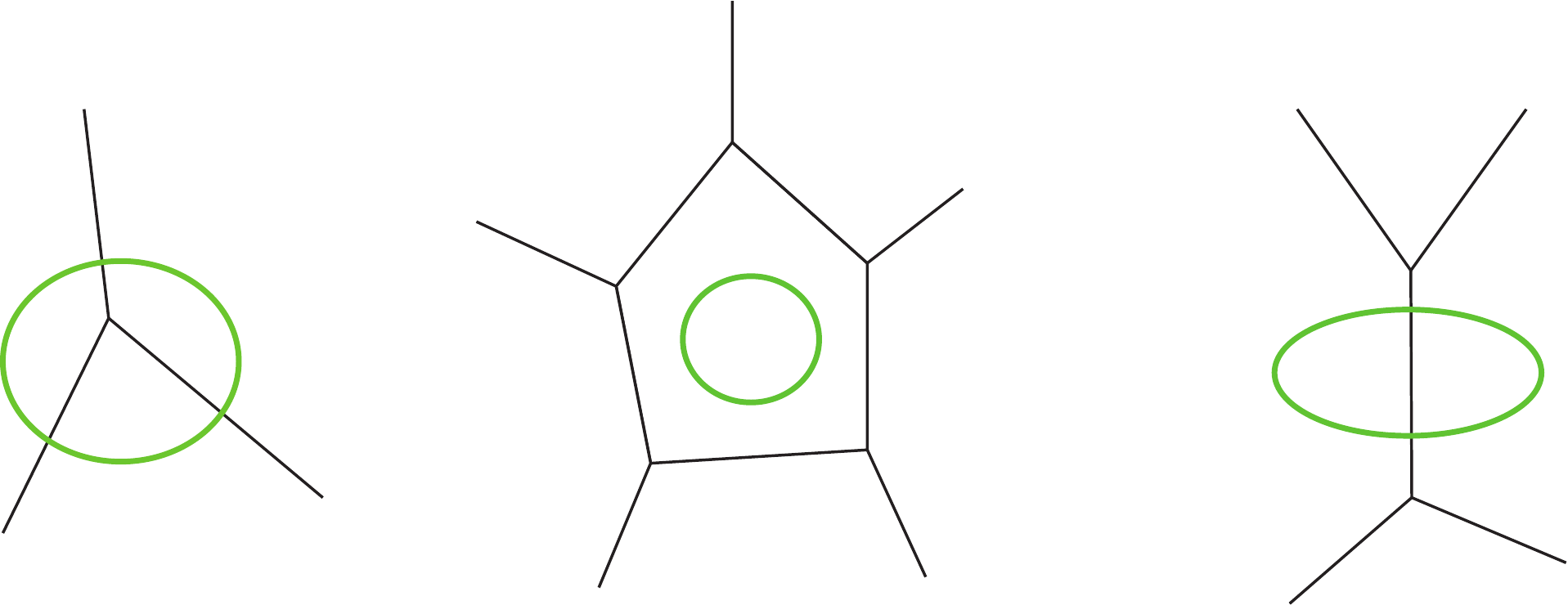} 
   \caption{What 2-orbifolds (in green) that are compressible look like}
   \label{fig:example}
\end{figure}

Figure 5 illustrates, in detail, what Definition 5 means. It shows the difference between what is compressible, what a compression is and what is not compressible. 

\begin{mydef} \emph{\textbf{(Separating Triangles)}}

A Separating Triangle is a prismatic 3-circuit. In an abstract polyhedron, if a Separating Triangle is encountered then the abstract polyhedron can be cut along the Separating Triangle. This will generate two separate polyhedron, both of which tile $\mathbb{H}^3$. 
\end{mydef}

\begin{figure}[htbp] 
   \centering
   \includegraphics[width=3.5in]{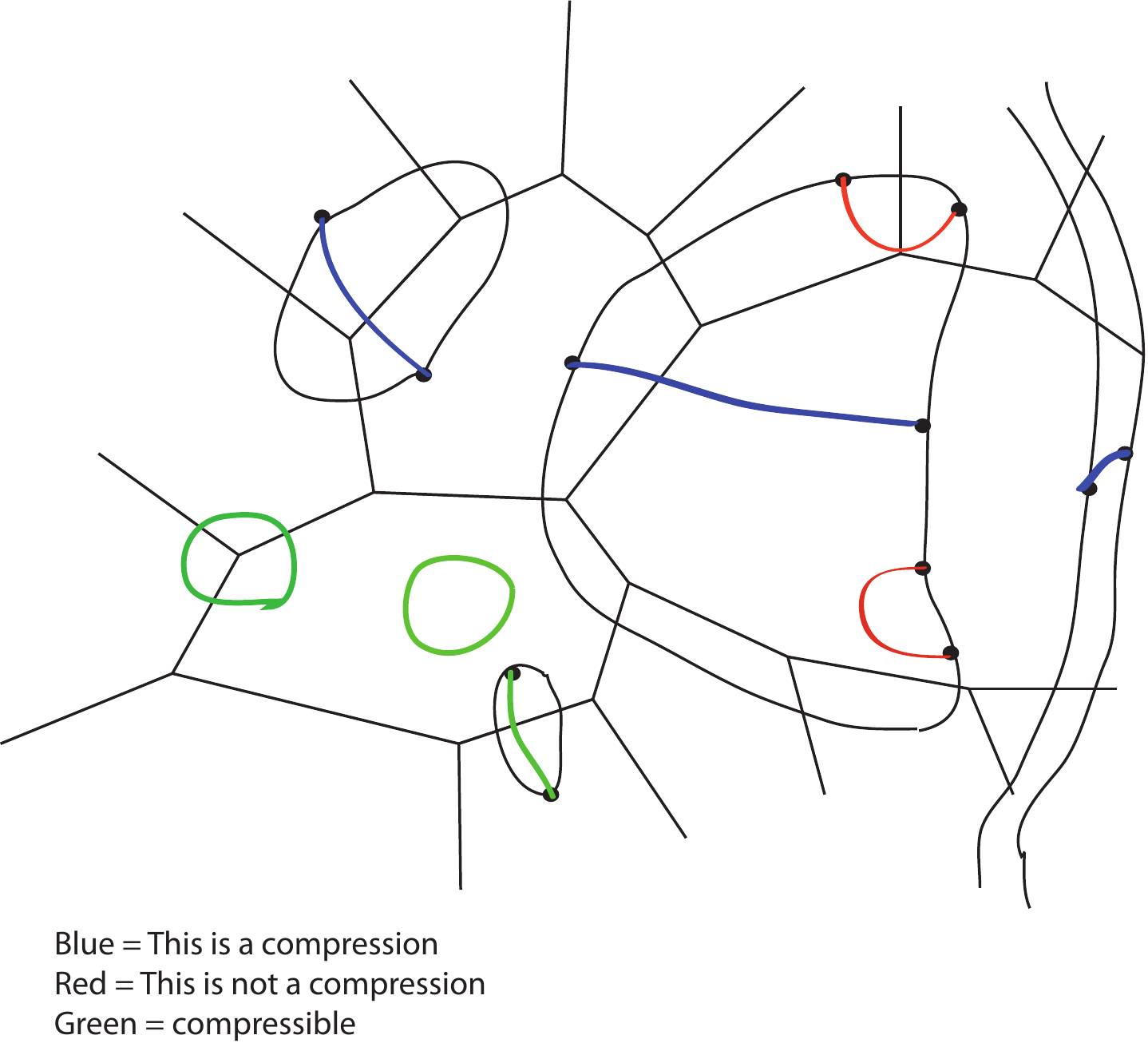} 
   \caption{The difference between compressible and a compression}
   \label{fig:example}
\end{figure}

\begin{mydef}
A polyhedron is called $\emph{large}$ if it contains an incompressible 2-orbifold and has no separating triangles. If the polyhedron does not contain an incompressible 2-orbifold and no separating triangles then it is considered to be small. 
\end{mydef}
 
A polyhedron that is considered to be small is a Tetrahedron. Large polyhedra are often called $\emph{Haken}$. 
The calculation of polyhedral volumes is complicated. For instance, the volume of an ideal tetrahedron (that is, a tetrahedron with all of whose vertices are at infinity) with dihedral angles $\alpha$, $\beta$, and $\gamma$ is $\cyrrm{D}(\alpha) + \cyrrm{D}(\beta) +\cyrrm{D}(\gamma)$ where $\cyrrm{d}(\theta) = - \int_0^\theta \ log | 2\sin(u) | du$ is the Lobachevsky function \cite{mil12}. 
Calculating the volume of a hyperbolic polyhedron is very difficult and requires a lot of time. Luckily, there is a computer program called Orb that does the calculations for us. Orb takes an abstract polyhedron with edge labels as input, and computes the volume of the polyhedron, if it is hyperbolic. Note: the number that is returns is double that of the actual volume. Orb will only measure the volume of an abstract polyhedron if it satisfies Andreev's Theorem. In the case that the labeled polyhedron does satisfy all of Andreev's conditions, Orb will return the volume with the word \textquotedblleft {Geometric}". When geometric is returned then that means that the volume is accurate. Below in Figure 4, an abstract polyhedron is shown with its approximate volume.

\begin{figure}[htbp] 
   \centering
   \includegraphics[width=5in]{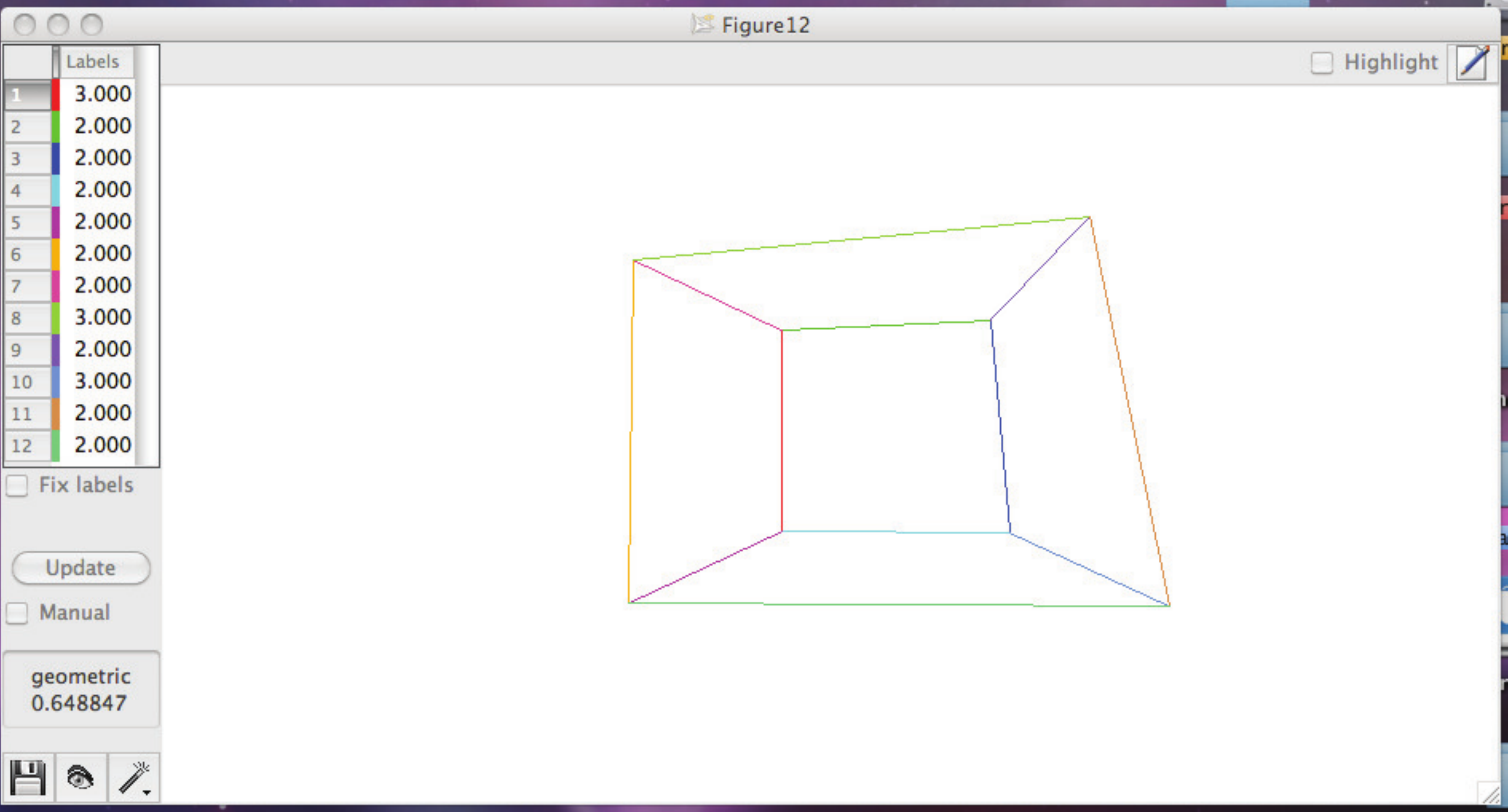} 
   \caption{Screen Shot of Orb with the Abstract polyhedron of a\textquotedblleft {cube} \textquotedblright} 
   \label{fig:example}
\end{figure}

\begin{myfact} 
Decreasing the angles of a hyperbolic polyhedra increases its \emph{volume}.
\end{myfact}

This fact can be easily witnessed by labeling the same polyhedron differently in Orb. The volume decreases as the label at any edge is increased. This is a consequence of Schl\"aflies formula.

Unlike Euclidean Geometry, in which is very easy to find volumes of polyhedra, hyperbolic Geometry poses quite the challenge, which is why Orb is so helpful. Orb relies on the Schl\"{a}flies equation in order to calculate the volume of polyhedra.

\begin{mydef} Schl\"{a}flies Formula for Polyhedra
\\
\\
Let $P_t$ be a smooth one-parameter family of polyhedra in $\mathbb{ H }^n$. Then the derivatives of the volume $V_t$ of $P_t$ satisfies the equation 
\begin{equation} 
-(n-1) \frac{dV}{dt} = \sum_{F}  \ V_{n-2} (F) \frac{d\theta_F}{dt}
\end{equation}
Where the sum is over all $(n-2)$-dimensional faces of $P_t$, and $V_{n-2}$ denotes the $(n - 2)$-dimensional volume and $\theta_t$ denotes the dihedral angle of $F$. \cite{CHK} 
\end{mydef}

In the two dimensional case, the equation is very simple.

\begin{equation}
-\frac{dA_t}{dt} = \sum \frac{d\theta_V}{d_t}
\end{equation}
where $\theta_t$ is the angle at the vertex and $A$ is the area.

\section{Investigations of Volumes}

I took the following as facts in my investigations:

\begin{enumerate}
\item A small, three dimensional Coxeter polyhedron has no separating triangles and no incompressible 2-orbifolds, except for those coming from triangular faces.
\item The larger the angle, i.e. the smaller the edge label, the higher the volume. 
\item An ideal vertex is created if the sum of the edges coming together in a vertex are exactly $\pi$.  
\item Every abstract polyhedron can be divided along its separating triangles into a polyhedra with no separating triangles. 
\item Hyperbolic Coxeter polyhedra are convex. Convex means that for any two points $P$ and $Q$ in the object, the shortest path from $P$ to $Q$ is totally contained in the object. 
\end{enumerate}

I started off by drawing abstract polyhedron and labeling them s0 that they satisfied all of Andreev's conditions. I then used Orb to calculate the volume. I noticed fairly early on that the higher the label of the edge, the higher the volume would be. Changing an edge label from a 5 to a 3 made a big difference, for instance. To see a planar graph of an object and know that is large was a little challenge, because sometimes it turned out that the abstract polyhedron that I drew had multiple separating triangles and once those were separated, the abstract polyhedron that I was left with was, in fact, small. Also, allowing only one vertex to go to infinity would increase the volume greatly. Chipping away at the faces and edges resulted in the smallest, large, Coxeter polyhedron, which turned out to have a planer graph of a cube. The cube that I had found already had a name, and it was referred to as  the \underline {Lambert cube}. The abstract polyhedron of the Lambert cube is shown in Figure 4.  

The smallest volume that could be determined for the Lambert cube was 0.648847. 

\subsection{Observations}
Before arriving at this conclusion, a couple of observations had to be made. First is that the Lambert cube, or at that time the Cube, contained three distinct prismatic 4-circuits. This is good, because that means that I had successfully removed all prismatic 3-circuits or separating triangles. This meant that the abstract polyhedron that I ended up with was in fact large.
\begin{figure}[htbp] 
   \centering
   \includegraphics[width=3in]{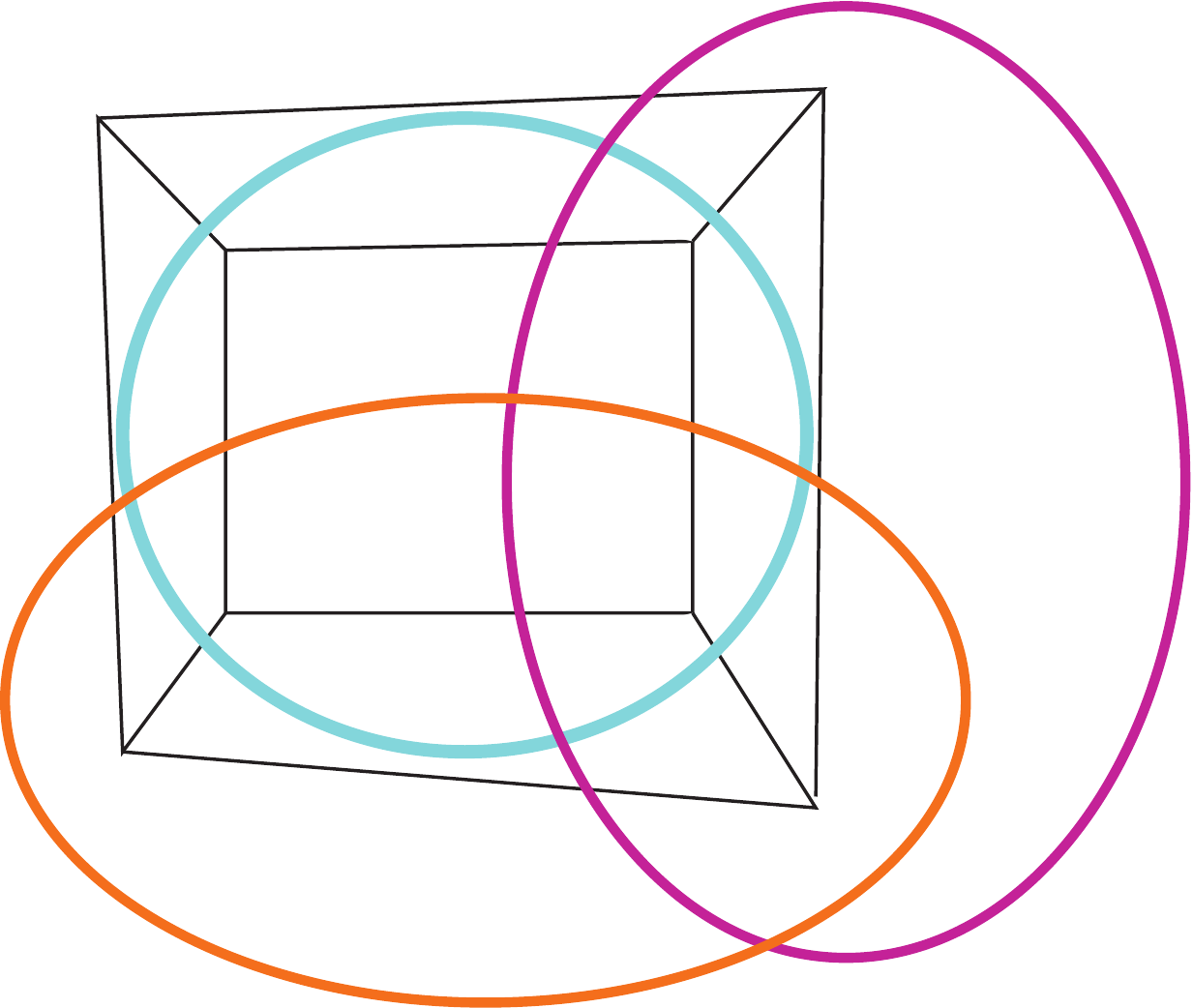} 
   \caption{The three prismatic 4-circuits on the Lambert Cube}
   \label{fig:example}
\end{figure}
\\
\\
\indent Due to Andreev's 4th condtion, the edges in the prismatic 4-circuits had to add up to less than $2\pi$. This means that at least one edge had to be labeled with a 3 or greater in each prismatic 4-circuit; in total there had to be at least three edges labeled with a 3 or greater. The placement of the 3s turned out to be also very important. If any of the 3s were placed adjacent to each other then the volume would increase. As already mentioned, if all of the 3s were surrounding one vertex, this created an ideal vertex, then the volume increased significantly.
\begin{figure}[htbp] 
   \centering
   \includegraphics[width=5in]{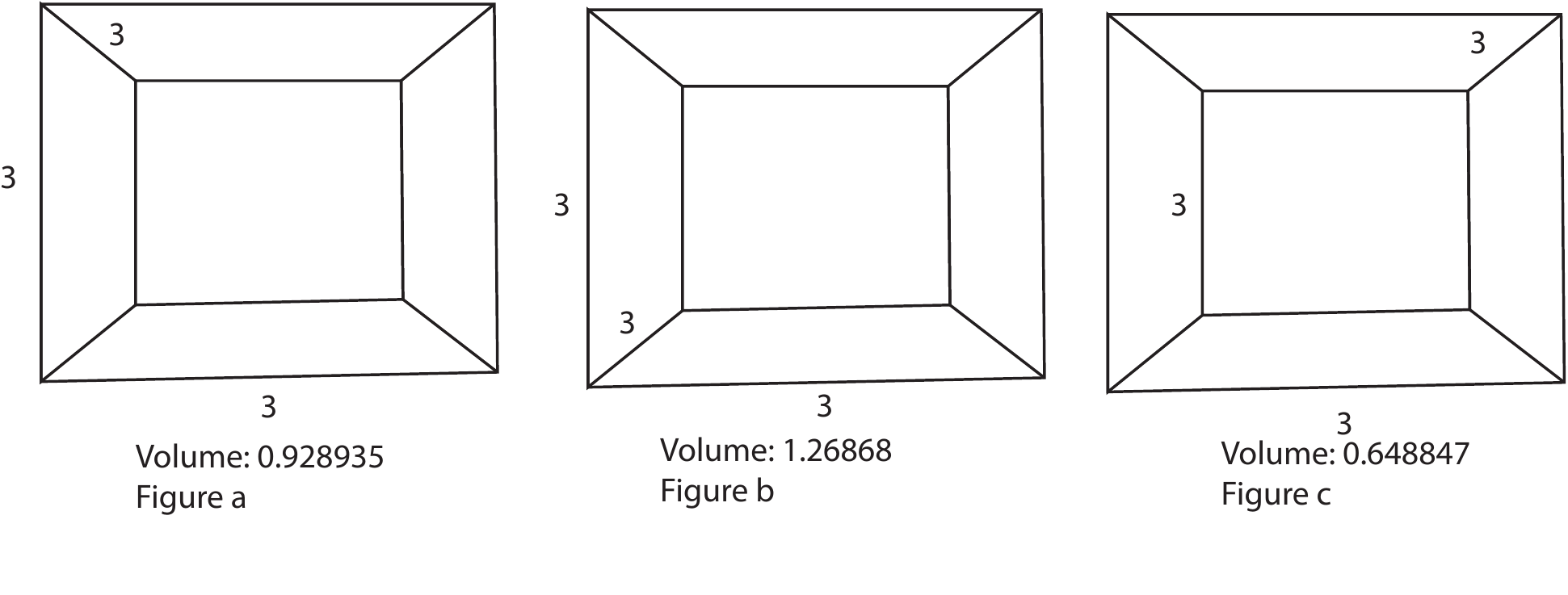} 
   \caption{Lambert Cube with different edge labelings and volumes. Edges not labeled have value 2. Figure b has one vertex at infinity.}
   \label{fig:example}
\end{figure}
\\
\\
\indent I reached the conclusion that I only needed 2s and 3s to enumerate the abstract polyhedron so that it would satisfy all of Andreev's conditions. There were multiple ways to place the three 3s without them being adjacent to each other. It turns out that it does not matter where the three 3s were placed so long as they were not adjacent. The volume would remain the same. I attribute this results to symmetry. 
 \\
\\
We summarize the following facts and observations:
\begin{enumerate}
\item Every polyhedron has a smallest possible volume. 
\item When a polyhedron is labeled with only the minimal amount of large numbers, the volume should be the smallest possible volume for that polyhedron.  
\item The smallest volume is obtained by moving the larger numbers as far away as possible from each other. This means that if there are no large numbers adjacent to each other then the volume should be the smallest. 
\item The amount of prismatic 3-circuits and 4-circuits which are distinct will determine the amount numbers larger than 3 are needed in the polyhedron to make it the smallest, Haken polyhedron. 
\end{enumerate}

Next, I started investigating Pyramids with one vertex at infinity. An abstract polyhedron of a Pyramid is illustrated in Figure 5 below. The edge which are going off to infinity all get labeled 2. 

\begin{figure}[htbp] 
   \centering
   \includegraphics[width=3.5in]{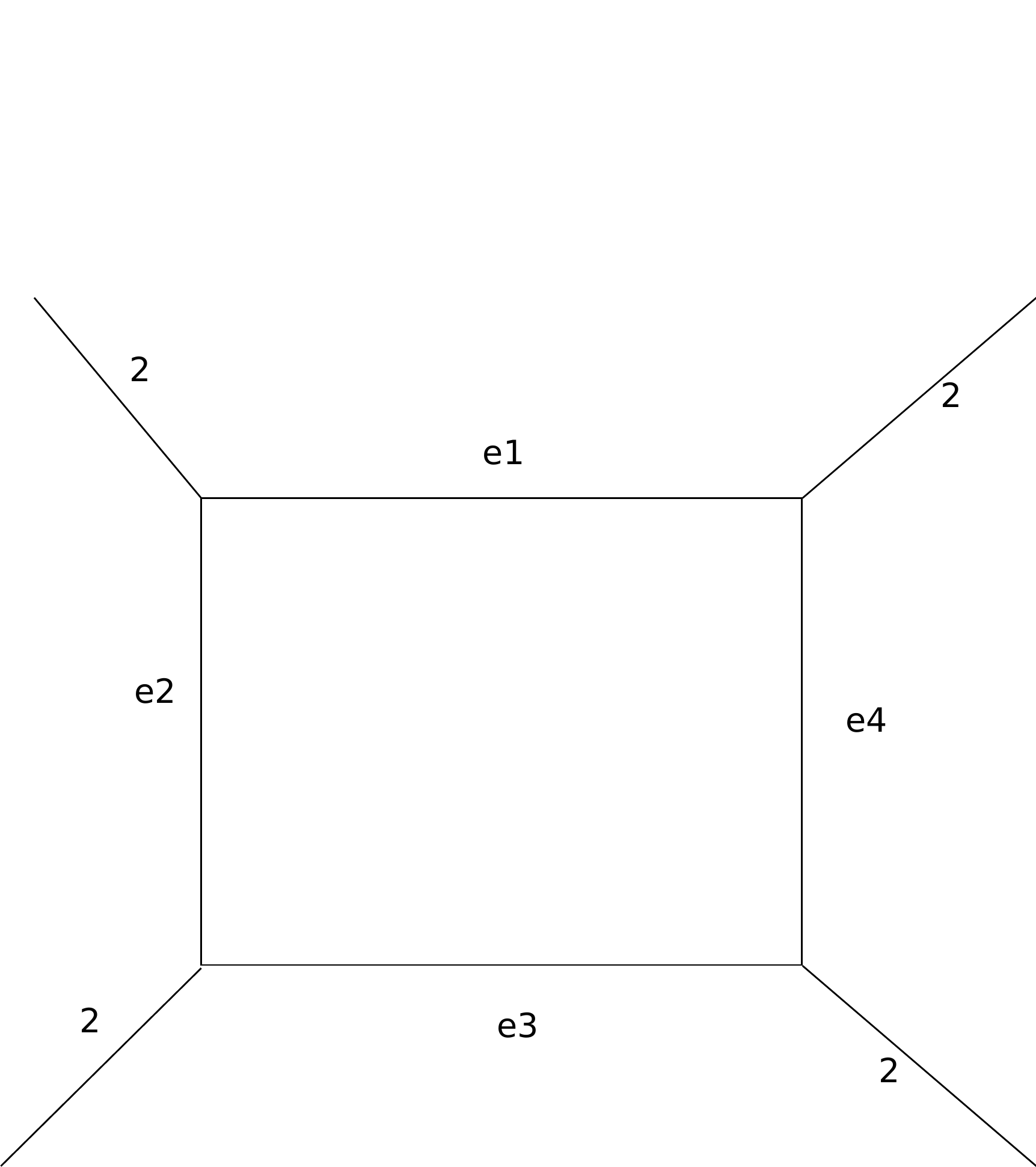} 
   \caption{Edge labeling of Pyramid with at least one vertex at Infinity }
   \label{fig:example}
\end{figure}

The remaining edges, $e_1, e_2, e_3, e_4$ are enumerated as shown in the table below. These are all possible combinations for which the corresponding polyhedron will which would tile $\mathbb{H}^3$ up to isometries. No other vertex can be at infinity. 
\begin{center}
\begin{tabular}{|l | r | r | r| }
\hline
$e_{1}$ 	& $e_{2}$	& $e_{3}$	& $e_{4}$ 	\\ \hline
2		& 2		&  3		& (3,4,5,6)	\\ \hline
2		& 2 		& 4		& 4		\\ \hline
2		& 3		& (3,4,5) 	& 3		\\ \hline
3		& (3,4,5)	& 3		& 3		\\ \hline
3		& (3,4,5)	& 3		& (3,4,5)	\\ \hline

\end{tabular} 
\end{center}
\vspace{ 1 cm }

\pagebreak

\section{Notes and Thanks}

\indent I would like to thank Michael Fry who really helped and worked hard on this project with me.\\ 
\indent I would also like to thank my advisor Professor Rafalski for all the hard work he put in this project for and with me as well as dragging me to Mathematics conferences. I also would like to thank my advisor again for all the support he has provided to me this past year and hopefully, that support will carry on after my undergraduate carrer. \\
\indent And last but not least, I would like to thank the entire Mathematics Department for always having an open and welcoming attitude towards me and my enthusiasm for the subject.

\end{document}